\documentstyle[12pt]{article}

\begin{document}
\begin{center}
{\Large {\bf Planscherel Measure on $E_q(2)$ }}
\end{center}

\vspace{1cm}

\begin{flushleft}
H. Ahmedov$^1$  and I. H. Duru$^{2,1}$

\vspace{.5cm}

{\small
1. TUBITAK, Feza G\"ursey Institute,  P.O. Box 6, 81220 Cengelkoy, 
Istanbul, Turkey 
\footnote{E--mail : hagi@fge1.gursey.gov.tr and duru@fge1.gursey.gov.tr}.

2. Trakya University, Mathematics Department, P.O. Box 126, 
Edirne, Turkey.}
\end{flushleft}

\vspace{2cm}

\begin{center}
{\bf Abstract}
\end{center}

Following the construction of the invariant integral and the scalar product
for the quantum Euclidean group $E_q(2)$, we obtained the full matrix
elements of its unitary irreducible representations from $SU_q(2)$ by
contraction and then derived the Planscherel measure.

\vspace{2cm}

\begin{center}
July 1998
\end{center}

\pagebreak

\vfill
\eject

\renewcommand{\theequation}{I.\arabic{equation}} \setcounter{equation}{0}

\section{Introduction}

If the space--time is not a differentiable manifold but a non--commutative
geometry there may be interesting new physical effects. Because of the non
existence of sufficiently developed differential calculus on the
non--commutative geometry to study physical effects one needs new
mathematical tools. Since on the other hand the quantum group spaces are
natural examples for the non--commutative geometries it is of interest to
develop an algorithm that can be employed in physical applications. In that
direction we have recently investigated the Green function on $SU_q(2)$ \cite
{1}. Purpose of the present work is to present some tools to be used in the
spectral analysis on group $E_q(2)$ or corresponding non--commutative plane $%
E_q(2)/U(1)$.

In the following section we first review the known invariant integral on $%
SU_q(2)$ and then present the similar construction for $E_q(2)$.

In Section III the scalar product on $E_q(2)$ is introduced.

Section IV is devoted to the derivation of the matrix elements of $E_q(2)$
irreducible representations from $SU_q(2)$ by contraction.

In Section V we obtain the Planscherel measure on $E_q(2)$.

Most of the known formulae as well as methods which we employ are presented
in the Appendix.

\newpage

\renewcommand{\theequation}{II.\arabic{equation}} \setcounter{equation}{0}

\section{Invariant Integral on $E_q(2)$}

Let us start by reviewing the construction of the invariant integral for $%
SU_q(2)$ which will guide in $E_q(2)$ case.

\vspace{.5cm}

\begin{flushleft}
{\bf Rewiev of the invariant integral on $SU_q(2)$}
\end{flushleft}

\vspace{.5cm}

The coordinate functions $x$, $u$, $x^*$, $u^*$ generating the $*$--Hopf
algebra $A(SU_q(2))$ satisfy the following relations \cite{2} 
\begin{eqnarray}\label{com}
u(x)=qx u, \ \  x^*u=qux^*,  \ \  uu^*=u^*u, \nonumber \\  
xx^* + uu^*={\bf 1}, \ \ x^*x+ q^2uu^*={\bf 1}. \ \ \ 
\end{eqnarray}
In $q\rightarrow 1$ limit the above relations define the three dimensional
sphere $S^3$ which is the topological manifold of $SU(2)$. In $q\neq 1$ case
we consider irreducible $*$--representation $\pi $ of the associative
algebra $A(SU_q(2))$ in the Hilbert space ${\cal L}^2(Z)$ with the
orthonormal basis $\{\mid n\rangle \}_{n\in Z}$ \cite{3} 
\begin{eqnarray}\label{rep} 
\pi (x)\mid n\rangle=(1-q^{2n})^{1/2}\mid n-1\rangle, \ \ \ \ \ \
\pi (u)\mid n\rangle=q^n\mid n\rangle, \nonumber \\
\pi (u^*)\mid n\rangle=q^{n}\mid n\rangle, \ \ \ \ \ \
\pi (x^*)\mid n\rangle=(1-q^{2n+2})^{1/2}\mid n+1\rangle
\end{eqnarray}
for $0<q<1$. We associate a vector $v\in {\cal L}^2(Z)$ to each point of the
``topological space" $SU_q(2)$. Then 
\begin{equation}
f_v=\langle v\mid \pi(f) \mid v \rangle 
\end{equation}
gives the value of any function $f\in A(SU_q(2))$ at this point. Thus by
topological space of the quantum group $SU_q(2)$ we mean the carrier space $%
{\cal L}^2(Z)$ of the $*$--representation of $A(SU_q(2))$. By the subspace $%
X_q$ of $SU_q(2)$ we mean the subspace $H^\prime$ of the Hilbert space $%
{\cal L}^2(Z)$. The empty set in $SU_q(2)$ is the empty set in ${\cal L}^2(Z)
$. Two subspaces $X_q$ and $X_q^\prime$ of $SU_q(2)$ are said to have zero
intersection ($X_q\bigcap X^\prime_q =\emptyset$ ) if the corresponding
subspaces $H$ and $H^\prime$ of the Hilbert space ${\cal L}^2(Z)$ has zero
intersection. The union and intersection of subspaces $X_q$, $%
X_q^\prime\subset SU_q(2)$ are understood as the union and intersection of
corresponding subspaces $H$, $H^\prime\subset {\cal L}^2(Z)$.

Consider the linear map $\mu : SU_q(2)\rightarrow [0, \infty )$ defined by 
\begin{equation}
\label{lin}\mu (X_q)=(1-q^2) \sum_{n\in J} \langle n\mid \pi (\mu_S )\mid n
\rangle, 
\end{equation}
where $J\subset Z$, such that the vectors $\mid n \rangle$, $n\in J$ span
the subspace $H$ of ${\cal L}^2(Z)$ corresponding to the subspace $X_q$ of $%
SU_q(2)$, $\mu_S\in A(SU_q(2))$. Since the left hand side of the above
expression is positive definite the operator $\pi(\mu_S)$ must be
self--adjoint and positive definite in ${\cal L}^2(Z_0)$. To make the linear
map ({\ref{lin}) a measure on $SU_q(2)$ we have to impose the additivity
condition 
\begin{equation}
\mu (\bigcup_j X_{qj}) =\sum_j \mu (X_{qj} ) 
\end{equation}
for the disjoint subspaces $X_{qj}$ of $SU_q(2)$. Inspecting (\ref{rep}) and
(\ref{lin}) we conclude that $\mu_S$ is a polynomial of $\xi =uu^* $. The
measure on $SU_q(2)$ is then given by 
\begin{equation}
\mu (X_q) =(1-q^2) \sum_{n\in J} \langle n\mid \pi (\mu_S(\xi ))\mid n
\rangle = (1-q^2) \sum_{n\in J} \mu_S(q^{2n}). 
\end{equation}
}

By means of the measure $\mu$ on the quantum group $SU_q(2)$ we introduce
the linear functional $\psi : A(SU_q(2))\rightarrow C$ as 
\begin{equation}
\label{linf}\psi (f) =(1-q^2) \sum_{n=0}^\infty \langle n\mid \pi
(f\mu_S)\mid n\rangle. 
\end{equation}
If $f$ is the function of $\xi$ only we can rewrite (\ref{linf}) as 
\begin{equation}
\label{linf1}\psi (f) = \int_{0}^1 f(\xi ) d_{q^2}\mu_S(\xi ). 
\end{equation}
On the other hand the invariant integral on $SU_q(2)$ is \cite{3} 
\begin{equation}
\psi (f) = \int_{0}^1 f(\xi ) d_{q^2}\xi . 
\end{equation}
Comparing to (\ref{linf1}) we see that if 
\begin{equation}
\label{eq: mej8}\mu_S(\xi )=\xi , 
\end{equation}
the linear functional (\ref{linf}) defines the invariant integral on $SU_q(2)
$. Hence, the invariant measure and invariant integral on the quantum group $%
SU_q(2)$ are given by 
\begin{equation}
\label{invm}\mu (X_q) = (1-q^2) \sum_{n\in J} q^{2n} 
\end{equation}
and 
\begin{equation}
\label{invin}\psi (f) =(1-q^2) \sum_{n=0}^\infty \langle n\mid \pi(f\xi
)\mid n\rangle 
\end{equation}
respectively. Here $J\subset Z$, such that the vectors $\mid n \rangle$, $%
n\in J$ span the subspace $H$ of ${\cal L}^2(Z)$ corresponding to the
subspace $X_q$ of $SU_q(2)$.

\vspace{1cm}

\begin{flushleft}
{\bf Invariant integral on $E_q(2)$ }
\end{flushleft}

\vspace{.5cm}

In fashion parallel to $SU_q(2)$ we can define the invariant measure on $%
E_q(2)$. The irreducible $*$--representation $\pi $ of the algebra of
polynomials on $E_q(2)$ is constructed in the Hilbert space ${\cal L}^2(S)$
of square integrable functions on the circle $S$. In the orthonormal basis $%
\mid j\rangle =\frac 1{\sqrt{2\pi }}e^{ij\psi }$, $-\infty <j<\infty $ we
have 
\begin{equation}
\label{repe1}\pi (z)\mid j\rangle =q^{-j}\mid j-1\rangle ,\ \ \ \pi
(z^{*})\mid j\rangle =q^{-j-1}\mid j+1\rangle ,\ \ \ \pi (\delta )\mid
j\rangle =\mid j-2\rangle .
\end{equation}
Thus as topological space the quantum Euclidean group $E_q(2)$ is equivalent
to the Hilbert space ${\cal L}^2(S)$. The measure on $E_q(2)$ is given by 
\begin{equation}
\mu (X_q)=(1-q^2)\sum_{j\in J}\langle j\mid \pi (\mu _E)\mid j\rangle ,
\end{equation}
where $\mu _E$ is polynomial of coordinate functions $z$, $z^{*}$, $\delta
^{\pm 1}$; and $J\subset Z$ such that the vectors $\mid j\rangle $, $j\in J$%
, span a basis in the subspace $H\subset {\cal L}^2(S)$ corresponding to the
subspace $X_q\subset E_q(2)$. Since $E_q(2)$ can be obtained from $SU_q(2)$
by contraction we define $\mu _E$ from $\mu _S$ of (\ref{invm}) as 
\begin{equation}
\mu _E=\lim _{r\rightarrow \infty }(r^2\mu _S(\frac{u_0u_0^{*}}{r^2}%
))=zz^{*}=\rho ^2.
\end{equation}
Thus the invariant measure on $E_q(2)$ is given by 
\begin{equation}
\mu (X_q)=(1-q^2)\sum_{j\in J}\langle j\mid \pi (\rho ^2)\mid j\rangle
=N\sum_{j\in J}q^{-2j}.
\end{equation}
By the virtue of the invariant measure we define the invariant integral on $%
E_q(2)$ as 
\begin{equation}
\psi (f)=(1-q^2)\sum_{j=-\infty }^\infty \langle j\mid \pi (f\rho ^2)\mid
j\rangle ,
\end{equation}
provided that the left--hand side is finite. If $f(g)=f(\rho )$ the above
expression can be rewritten by means of $q$--integral as 
\begin{equation}
\psi (f)=\int_{-\infty }^\infty f(\rho )d_{q^2}(\rho ^2).
\end{equation}
The invariant measure and invariant integral on the quantum group $E_q(2)$
are then given by 
\begin{equation}
\label{invme}\mu (X_q)=(1-q^2)\sum_{j\in J}q^{-2j}
\end{equation}
and 
\begin{equation}
\label{invine}\psi (f)=(1-q^2)\sum_{j=-\infty }^\infty \langle j\mid \pi
(f\rho ^2)\mid j\rangle 
\end{equation}
respectively. The subset $J\subset Z$ is defined such that the vectors $\mid
j\rangle $, $j\in J$ span the subspace $H\subset {\cal L}^2(S)$
corresponding to the subspace $X_q$ of $E_q(2)$.

\vspace{1cm}

\renewcommand{\theequation}{III.\arabic{equation}} \setcounter{equation}{0}

\section{Scalar Product on $E_q(2)$}

Let $\Phi (E_q(2))$ be the set of analytic functions on $E_q(2)$ such that 
\begin{equation}
\label{plan1}\psi (f f^*)<\infty , 
\end{equation}
where $\psi$ is the invariant integral (\ref{invine}) on $E_q(2)$.

Recall that the homomorphism 
\begin{equation}
\phi (z)=0, \ \ \ \phi (z^*)=0, \ \ \ \phi (\delta)=t 
\end{equation}
defines the quantum subgroup $U(1)\subset E_q(2)$. We have then the
decomposition 
\begin{equation}
\label{plan2}\Phi(E_q(2))=\sum_{ij\in Z}\oplus \Phi[i,j] 
\end{equation}
where 
\begin{equation}
\Phi[i,j]=\{ f\in \Phi(E_q(2)) : L (f) = t^i\otimes f ; \ \ R (f) = f\otimes
t^j \} 
\end{equation}
and 
\begin{equation}
L = (\phi \otimes id)\circ \Delta ,\ \ \ \ R = (id\otimes \phi )\circ\Delta
. 
\end{equation}
The subspace $\Phi[i,j]\subset\Phi(E_q(2))$ consists of the elements of the
following form 
\begin{equation}
\label{plan3}f_{ij} (g) = \delta^jz^{i-j} f_{ij}(\rho^2), \ \ \ \ \ for \ \
i\geq j 
\end{equation}
and 
\begin{equation}
\label{plan4}f_{ij} (g) = \delta^j(z^*)^{j-i} f_{ij}(\rho^2), \ \ \ \ \ for
\ \ i\leq j, 
\end{equation}
where $\rho^2 =zz^*$.

By means of the invariant integral we introduce in $\Phi (E_q(2))$ the
bilinear forms 
\begin{equation}
\label{plan6}(f , f^\prime )_L =\psi (f^* f^\prime ) 
\end{equation}
and 
\begin{equation}
\label{plan7}(f^\prime , f )_R =\psi (f^\prime f^*) 
\end{equation}
related to each other as 
\begin{equation}
\label{plan8}( f , f^\prime )_L=(\tau(f^\prime ), f)_R, 
\end{equation}
where $\tau $ is the automorphism in $E_q(2)$ defined as 
\begin{equation}
\tau (z)=qz, \ \ \tau (z^*)=q^{-2}z^*, \ \ \tau (\delta )=q^{-4}\delta 
\end{equation}
By the virtue of (\ref{plan3}) and (\ref{plan4}) we get 
\begin{equation}
\label{top4}\tau (f)=q^{-2(i+j)}f 
\end{equation}
for $f\in \Phi[i,j]$. To prove the identity (\ref{plan8}) we use the
following representation for the invariant integral 
\begin{equation}
\label{plan10}\psi (f)=(1-q^2)Tr (f\rho^2), 
\end{equation}
which is the result of (\ref{invine}). In this representation the equality (%
\ref{plan8}) reads 
\begin{equation}
Tr( f^* f^\prime \rho^2)=Tr(\tau(f^\prime )f^*\rho^2). 
\end{equation}
Due to the decomposition (\ref{plan2}) it is enough to verify that (\ref
{plan10}) is valid for $f\in \Phi [i,j]$ and $f^\prime\in\Phi [i^\prime,
j^\prime]$. The latter can easily be verified by using (\ref{plan3}), (\ref
{plan4}), (\ref{top4}) and (B.7).

The representation ({\ref{plan10}) of the invariant integral allows us to
prove that the bilinear forms (\ref{plan6}) and (\ref{plan7}) are scalar
products in $\Phi (E_q(2))$. For that purpose we put $F=f\rho $ and $%
F^{\prime }=f^{\prime }\rho $ in (\ref{plan6}) and using (\ref{plan10}) we
get 
\begin{equation}
(f,f^{\prime })_L=(1-q^2)Tr(F^{*}F^{\prime }).
\end{equation}
The left hand side of the above equality defines the scalar product in the
space of Hilbert--Schmidt type operators. Thus the bilinear form (\ref{plan6}%
) defines the scalar product in $\Phi (E_q(2))$. In a similar fashion one
can show that bilinear form (\ref{plan7}) is also scalar product. }

From (\ref{plan3}), (\ref{plan4}), (\ref{repe1}) and (\ref{invine}) we have 
\begin{equation}
(f_1 , f_2 )_{L, R} =0, 
\end{equation}
for $f_1\in \Phi [i,j]$, $f_1\in \Phi [i^{\prime},j^{\prime}]$ such that $%
(i,j)\neq (i^{\prime},j^{\prime})$. Thus the decomposition (III.3) is
orthogonal with respect to the scalar products (\ref{plan6}) and (\ref{plan7}%
).

\vspace{1cm}

\renewcommand{\theequation}{IV.\arabic{equation}} \setcounter{equation}{0}

\section{Matrix Elements of the Irreducible Representations of $E_q(2)$ from 
$SU_q(2)$ by Contraction}

Unitary representations of $E_q(2)$ are previously studied \cite{5}. We now
show that the matrix elements of the unitary irreducible representations of $%
E_q(2)$ can be obtained from the ones of $SU_q(2)$ by contraction. This
example shows that there exists q--analog of contraction procedure of
classical groups representations which was investigated in \cite{6}.

In Appendix C we review the derivation of the matrix elements of the unitary
irreducible representations of $E(2)$ from $SU(2)$ by the following
contraction procedure 
\begin{equation}
t_{ij}^p (\phi , \rho , \zeta )= \lim _{l\rightarrow \infty }t_{ij}^l(\phi
,p\rho /l,\zeta-\phi ), 
\end{equation}
where one puts $\theta =p\rho /l$ in the decomposition (C.5) of $SU(2)$. The
quantum analog of the above procedure is 
\begin{equation}
\label{eq: c4}t_{ij}^p=\lim _{l\rightarrow \infty } t_{ij}^l(x_0, y_0, \frac{%
pv_0}{[l]}, \frac{pu_0}{[l]}), 
\end{equation}
where $[m]=\frac{q^m-q^{-m}}{q-q^{-1}}$ and $t^l_{ij}(x,y,v,u)$ are the
matrix elements of the unitary irreducible representation of $SU_q(2)$ (see 
\cite{vil} and references therein ) : 
\begin{equation}
t^l_{ij} =\lambda^l_{ij}x^{-i-j}v^{i-j}  _2\phi_1(q^{-2(l+j)},q^{2(j+l+1)};
q^{2(l+i-j)}\mid q^2, - q^2uv), 
\end{equation}
with $i+j\leq 0$, $j\leq i$ and 
\begin{equation}
\lambda^l_{ij} =q^{(l+i)(l-j)} \sqrt{\left[ 
\begin{array}{c}
l+j \\ 
j-i 
\end{array}
\right ]_{q^2} \left[ 
\begin{array}{c}
l-j \\ 
j -i 
\end{array}
\right ]_{q^2}}. 
\end{equation}
Here $_2\phi_1$ and $[\cdot]_q$ are the q--hypergeometric function and the
q--binomial coefficients respectively. We first calculate 
\begin{equation}
\lim _{l\rightarrow \infty } \lambda_{i,j}^lx_0^{-i-j} (\frac{(pv)}{[l]}%
))^{i-j} = \frac{q^{\frac{m^2-n^2}{2}+ \frac{j-i}{2}}}{[i-j]!}%
p^{i-j}x_0^{-i-j}v_0^{i-j} 
\end{equation}
and 
\begin{eqnarray}
\lim_{l\rightarrow\infty }
(\phi_{21}(q^{-2(l+j)},q^{2(l+1+j)};q^{2(1+i-j)}\mid q^2,
q^2 \frac{ (p^2v_0v_0^*) }{[l]})  ) = 
\nonumber \\
=[i-j]! \sum_{k=0}^\infty\frac{(-1)^k}{[k]![k+i-j]!}
(q^{-i-j}p^2zz^* )^k.
\end{eqnarray}  
We then combine them to arrive at 
\begin{equation}
\label{irep12}t^p_{ij}(g)=(iq^{-1/2})^{i-j}\delta ^{-j/2}(pz^*)^{i-j}{\cal J}%
_{i-j}(p^2zz^*)\delta ^{-j/2}, \ \ \ i\geq j. 
\end{equation}
For $i\leq j $ on the other hand one obtains 
\begin{equation}
\label{irep13}t^p_{ij}(g)=(-iq^{1/2})^{i-j}\delta ^{-j/2}{\cal J}%
_{j-i}(p^2zz^*)(pz)^{j-i}\delta ^{-j/2}. 
\end{equation}
In the above formulae ${\cal J}_j$ are the q--Bessel functions given by  
\begin{equation}
{\cal J}_j(x)=\sum_{k=0}^\infty \frac{(-1)^k} {[k]![k+j]!} (q^{-j}x)^k. 
\end{equation}

Before closing this section for the sake of completeness we like to present
the already known formulae for the right and left representations of the
quantum algebra $U_q(e(2))$ obtained from the representation of $U_q(su(2))$
by contraction \cite{dob} 
\begin{equation}
\label{eq: c6}{\cal R}(E_{\pm })t_{ij}^p=pt_{i\pm 1,j}^p,\ \ \ {\cal R}%
(k)t_{ij}^p=q^{-i}t_{ij}^p
\end{equation}
and 
\begin{equation}
\label{eq: c8}{\cal L}(E_{\pm })t_{ij}^p=pt_{i,j\mp 1}^p,\ \ \ {\cal L}%
(k)t_{ij}^p=q^{-j}t_{ij}^p.
\end{equation}

\vspace{1cm} \renewcommand{\theequation}{V.\arabic{equation}} 
\setcounter{equation}{0}

\section{Planscherel Measure on $E_q(2)$}

The comultiplication 
\begin{equation}
\label{eq: plan11}\Delta :\Phi (E_q(2))\rightarrow \Phi (E_q(2))\otimes \Phi
(E_q(2))
\end{equation}
defines the regular representations of the quantum group $E_q(2)$ in $\Phi
(E_q(2))$. Since the scalar product in $\Phi (E_q(2))$ is defined by means
of the invariant integral this representation is unitary. The linear space $%
\Phi (E_q(2))$ is common invariant dense domain for the set of linear
operators ${\cal R}(\phi )$, ${\cal L}(\phi )$, $\phi \in U_q(e(2))$. Since
the representation (\ref{eq: plan11}) of $E_q(2)$ is unitary the
representatives of ${\cal R}(E_{+}E_{-})$, ${\cal R}(H)$, ${\cal L}(H)$ will
be at least symmetric operators in $\Phi (E_q(2))$. From (\ref{eq: c6}), (%
\ref{eq: c8}) as well as from 
\begin{equation}
{\cal R}(E_{+}E_{-})t_{ij}^p=p^2t_{ij}^p
\end{equation}
we see that the matrix elements $t_{ij}^p$ are eigenfunctions of these
operators. Thus, the eigenfunctions $t_{ij}^p$ and $t_{i^{\prime }j^{\prime
}}^{p^{\prime }}$ corresponding to different eigenvalues are orthogonal 
\begin{equation}
(t_{ij}^p,t_{i^{\prime }j^{\prime }}^{p^{\prime }})_L=c_{ij}^l(p)\delta
_{ii^{\prime }}\delta _{jj^{\prime }}\delta (p-p^{\prime }),
\end{equation}
\begin{equation}
(t_{ij}^p,t_{i^{\prime }j^{\prime }}^{p^{\prime }})_R=c_{ij}^r(p)\delta
_{ii^{\prime }}\delta _{jj^{\prime }}\delta (p-p^{\prime }),
\end{equation}
where $c_{ij}^l(p)$ and $c_{ij}^r(p)$ are the normalization constants. Due
to the existence of the delta function the matrix elements of the unitary
irreducible representations of $E_q(2)$ in $\Phi (E_q(2))$ do not belong to $%
\Phi (E_q(2))$. This is natural feature of the non--compact groups.

Since $\tau (t_{ij}^p)=q^{2(i+j)} t_{ij}^p$ from (\ref{plan8}) and the above
orthogonality conditions we obtain 
\begin{equation}
\label{eq: plan14}c_{ij}^l(p)=q^{2(i+j)}c_{ij}^r(p). 
\end{equation}
From the explicit expressions (\ref{irep12}) and (\ref{irep13}) for the
matrix elements we have 
\begin{equation}
t_{ij}^{pp_0}=\beta (t_{ij}^p) 
\end{equation}
with $\beta $ being the automorphism in the quantum group $E_q(2)$ defined
as 
\begin{equation}
\beta (z)=p_0z,\ \ \ \beta (\delta )=\delta ,\ \ \ \beta (z^*)=p_0z^*;\ \ \
\ p_0\in (0,\ \infty ). 
\end{equation}
To preserve the unitarity of $E_q(2)$ representations we choose $p_0$ to be
real and positive. Using the explicit expressions for the Hermitian forms (%
\ref{plan6}) and (\ref{plan7}) we arrive at 
\begin{equation}
( t_{ij}^p, t_{i^{\prime }j^{\prime }}^{p^{\prime}}) _{R,L}= p_0^2(
t_{ij}^{pp_0}, t_{i^{\prime}j^{\prime }}^{p^{\prime }p_0})_{R,L} 
\end{equation}
which implies 
\begin{equation}
c_{ij}^{r,l}(p)=p_0c_{ij}^{r,l}(pp_0) 
\end{equation}
for any $p_0\in (0,\infty )$ and $i,j\in Z$. Thus the normalization
coefficients can be represented as 
\begin{equation}
c_{ij}^{r,l}(p)=\frac{1}{p}a_{ij}^{r,l}, 
\end{equation}
where the coefficients $a_{ij}^{r,l}$ do not depend on $p$.

From the unitarity of the left regular representation we have 
\begin{equation}
({\cal R}(E_\pm^n) t_{ij}^p, t_{i^\prime j^\prime }^{p^\prime})_R
=(t_{ij}^p, {\cal R}(E_\mp^n)t_{i^\prime j^\prime }^{p^\prime})_R. 
\end{equation}
The above relation implies 
\begin{equation}
a_{i,j}^r=a^r_{i \pm n,j} 
\end{equation}
for any $n\in Z_0$. Thus the coefficients $a_{i,j}^r$ do not depend on the
first subindex $i$. 
\begin{equation}
a_{i,j}^r=c^r_j 
\end{equation}
In a similar fashion using the unitarity of the right regular representation
we get 
\begin{equation}
a_{i,j}^l=c^l_i 
\end{equation}
which is independent of its second subindex. By the virtue of (\ref{eq:
plan14}) we have 
\begin{equation}
c_i^l=q^{2(i+j)}c_j^r, 
\end{equation}
which is solved by $c_i^l=cq^{2i}$ and $c_j^r=cq^{-2j}$ with $c=constant$.
Thus the matrix elements of the unitary irreducible representations satisfy
the following orthogonality conditions 
\begin{equation}
\label{eq: plan15}( t_{ij}^p, t_{i^{\prime }j^{\prime }}^{p^{\prime }})_R=%
\frac{cq^{-2j}}p\delta _{ii^{\prime }}\delta _{jj^{\prime }} \delta
(p-p^{\prime }), 
\end{equation}
and 
\begin{equation}
\label{eq: plan16}(t_{ij}^p\mid t_{i^{\prime }j^{\prime }}^{p^{\prime }}) _L=%
\frac{cq^{2i}}p\delta _{ii^{\prime }}\delta _{jj^{\prime }}\delta
(p-p^{\prime }), 
\end{equation}
which implies that the Planscherel measure on $E_q(2)$ is $p$.

For any function $f\in \Phi (E_q(2))$ which can be expressible as the linear
combination of the matrix elements $t_{ij}^p$ we have 
\begin{equation}
f=\frac 1c\int_0^\infty pdp\sum_{i,j=-\infty }^\infty q^{2j}\hat
f_{ij}^pt_{ij}^p,
\end{equation}
where 
\begin{equation}
\hat f_{ij}^p=(f,t_{ij}^p)_R.
\end{equation}

\newpage

\begin{center}
{\large {\bf Appendix}}
\end{center}

\renewcommand{\theequation}{A.\arabic{equation}} \setcounter{equation}{0}

{\bf A. Quantum Group $SU_q(2)$ and Algebra $U_q(su(2))$}

\vspace{.5cm}

The quantum group $SU_q(2)$ or the $*$--Hopf algebra $A(SU_q(2))$ is the
algebra of polynomials of the coordinate functions $x$, $x^*$, $u$ and $u^*$
satisfying the relations (\ref{com}) and the coalgebra operations 
\begin{equation}
\Delta (x)=x\otimes x -qu\otimes u^*, \ \ \ \Delta (u)=x\otimes u +u\otimes
x^*, 
\end{equation}
\begin{equation}
\varepsilon (x)=1, \ \ \ \varepsilon (u)=0 
\end{equation}
and the antipode 
\begin{equation}
S(x)=x^*, \ \ \ S(x^*)=x, \ \ \ S(u)=-qu, \ \ \ S(v)=-q^{-1}v. 
\end{equation}

The quantum algebra $U(su_q(2))$ is generated by the elements 
\begin{equation}
{\cal E_\pm}, \ k_\pm=q^{\pm H/4} 
\end{equation}
satisfying the relations 
\begin{equation}
[{\cal E}_{+},{\cal E}_{-}]= \frac{k^2-k^{-2}}{q-q^{-1}}, \ \ \ k{\cal E}%
_\pm = q^{\pm 1}{\cal E}_\pm k 
\end{equation}
involutions 
\begin{equation}
({\cal E}_\pm)^*={\cal E}_\mp , \ \ \ \ k^*=k 
\end{equation}
and co--algebra operations 
\begin{equation}
\Delta_U ({\cal E}_\pm )={\cal E}_\pm \otimes k +k^{-1}\otimes {\cal E}_\pm
, \ \ \ \Delta_U (k)=k\otimes k, 
\end{equation}
\begin{equation}
\varepsilon_U ({\cal E}_\pm )=0 , \ \ \ \ \varepsilon_U (k )=1 
\end{equation}
and antipode 
\begin{equation}
S_U ({\cal E}_\pm )=-q^{\pm 1}{\cal E}_\pm , \ \ \ \ S_U (k )=k^{-1}. 
\end{equation}
The quantum algebra $U_q(su(2))$ is in the non--degenerate duality with 
the quantum group $SU_q(2)$ 

\vspace{1cm} \renewcommand{\theequation}{B.\arabic{equation}} 
\setcounter{equation}{0}

{\bf B. $E_q(2)$ from $SU_q(2)$ by Contraction}

\vspace{.5cm}

Substituting 
\begin{equation}
\label{eq: contr1}x\rightarrow x_0, \ \ \ u\rightarrow \frac{1}{r}u_0, 
\end{equation}
in the formulas (II.1) and (A.1)$\rightarrow$(A.3) we get the relations 
\begin {eqnarray}
u_0u^*_0=u^*_0u_0, \ \  u_0x_0=qx_0u_0, \ \ \ qu_0x^*_0=x^*_0u_0, \nonumber \\ 
 x_0x^*_0+\frac{1}{r^2}u_0u^*_0=1,\ \ \ x^*_0x_0+\frac{q^2}{r^2}u_0u^*_0=1,
\end{eqnarray}  
the coalgebra operations 
\begin{equation}
\Delta x_0=x_0\otimes x_0 - \frac{q}{r^2}u_0\otimes u^*_0, \ \ \ \Delta
u_0=x_0\otimes u_0+ u_0\otimes x^*_0, 
\end{equation}
\begin{equation}
\varepsilon (x_0)=1,\ \ \ \varepsilon (u_0)=0 
\end{equation}
and the antipode 
\begin{equation}
S(x_0)=x^*_0,\ \ S(u_0)=-qu_0. 
\end{equation}
Taking $r\rightarrow \infty$ limit in the above formulae and putting 
\begin{equation}
\label{eq: contr2}z=iqx_0 u_0, \ \ \ \delta =x_0^2 
\end{equation}
we arrive at 
\begin{equation}
\label{eq: contr3}zz^*=q^{-2}z^*z,\ \ \ z\delta =q^2\delta z,\ \ \ z^*\delta
=q^2\delta z^*, 
\end{equation}
\begin{equation}
\label{eq: contr4}\delta ^{*}=\delta ^{-1}, 
\end{equation}
\begin{equation}
\label{eq: contr5}\Delta (z^*)=z^*\otimes 1+\delta^{-1} \otimes z^*, \ \ \
\Delta (z)=z\otimes 1+\delta \otimes z, \ \ \ \Delta (\delta )=\delta\otimes
\delta , 
\end{equation}
\begin{equation}
\label{eq: contr6}\varepsilon (\delta )=1,\ \ \ \varepsilon (z)=0, 
\end{equation}
and 
\begin{equation}
\label{eq: contr7}S(\delta )=\delta ^{-1},\ \ \ S(z)=-\delta ^{-1}z,\ \ \
S(z^*)=-\delta z^*. 
\end{equation}
The above relations define the quantum Euclidean group $E_q(2)$ \cite{cel}.

Due to the duality one obtains $U_q(e(2))$ from $U_q(su(2))$ also by
contraction. Substituting 
\begin{equation}
\label{eq: contr10}{\cal E}_\pm \rightarrow r E_\pm, \ \ \ k\rightarrow k 
\end{equation}
in (A.5)$\rightarrow$(A.9) and taking $r\rightarrow\infty$ limit one gets
the relations 
\begin{equation}
\label{eq: contr12}[ E_{+}, E_{-}]=0, \ \ \ k E_\pm =q^{\pm 1}E_\pm k, 
\end{equation}
the involution 
\begin{equation}
\label{eq: contr13} E_\pm^*=E_\mp, \ \ \ k^*=k, 
\end{equation}
the coalgebra operations 
\begin{equation}
\Delta_U ( E_\pm )= E_\pm \otimes k +k^{-1}\otimes E_\pm , \ \ \ \Delta_U
(k)=k\otimes k, 
\end{equation}
\begin{equation}
\varepsilon_U ( E_\pm )=0 , \ \ \ \ \varepsilon_U (k )=1 
\end{equation}
and the antipode 
\begin{equation}
S_U ( E_\pm )=-q^{\pm 1} E_\pm , \ \ \ \ S_U (k )=k^{-1}, 
\end{equation}

\vspace{1cm}

\renewcommand{\theequation}{C.\arabic{equation}} \setcounter{equation}{0}

{\bf C. Unitary Representations of $E(2)$ from $SU(2)$ by Contraction \cite
{9}}

\vspace{.5 cm}

Matrix elements $t^p_{kj}$, $p\in (0, \ \infty )$, $-\infty < k, \ j < \infty
$, of the unitary irreducible representation of $E(2)$ has the following
integral representation 
\begin{equation}
t_{kj}^p(g_E) = \frac{e^{-i (k\phi +j(\zeta -\phi))}}{2\pi }
\int_0^{2\pi}e^{ip\rho \cos \psi } e^{i(k-j)\psi }d\psi , 
\end{equation}
where 
\begin{equation}
g_E= \left ( 
\begin{array}{ccc}
\cos\zeta & -\sin\zeta & \rho\cos\phi \\ 
\sin\zeta & \cos\zeta & \rho\sin\phi \\ 
0 & 0 & 1 
\end{array}
\right ). 
\end{equation}
Matrix elements $t^l_{ij}$, $l\in\frac{1}{2}Z_0$, $-l \leq i,j\leq l$, of
the unitary irreducible representation of $SU(2)$ are given by 
\begin{equation}
t_{kj}^l(g_S)=e^{-i(k\phi +j\phi ^\prime ) }P_{kj}^l(\cos \theta ), 
\end{equation}
where $P_{kj}^l$ is the Jacobi polynomial which has the following integral
representation 
\begin{eqnarray}
\label{eq: c1} P_{kj}^l(\cos\theta )=
\frac{1}{2\pi}\sqrt{\frac{(l-j)!(l+j)!}{(l-k)!(l+k)!}} \int_0^{2\pi} d\psi
e^{ij\psi}
 (i\sin(\theta /2)e^{i\psi /2} 
+\cos(\theta /2)e^{-i\psi /2} )^{l-k} \nonumber \\  
(i\sin(\theta /2)e^{-i\psi /2} +\cos(\theta /2)e^{i\psi /2})^{l+k} \ \ \ \ \ 
\ \  \ \   \ \ \ \ 
\end{eqnarray}
and 
\begin{equation}
\label{eq: c0}g_S= \left( 
\begin{array}{cc}
e^{i(\phi +\phi^\prime) /2}\cos \theta /2 & ie^{i(\phi-\phi^\prime ) /2}\sin
\theta /2 \\ 
ie^{i(\phi^\prime -\phi) /2}\sin \theta /2 & e^{-i(\phi +\phi^\prime)
/2}\cos \theta /2 
\end{array}
\right) 
\end{equation}
Putting $\theta =p\rho /l$ in (\ref{eq: c1}) for $l>>1$ we have 
\begin{equation}
P_{kj}^l(\cos(p\rho /l))=\frac{1}{2\pi} \int_0^{2\pi} d\psi
e^{i(j-k)\psi}d\psi (1 +\frac{ip\rho }{2l}e^{-i\psi /2} )^{l-k} (1 +\frac{%
ip\rho }{2l}e^{i\phi /2} )^{l+k}. 
\end{equation}
By the virtue of $\lim_{l\rightarrow \infty} (1+x/l)^l = e^x$ we obtain 
\begin{equation}
\lim_{l\rightarrow \infty}P_{kj}^l(\cos(p\rho /l))= \frac{1}{2\pi}
\int_0^{2\pi}e^{ip\rho \cos \psi} e^{i(k-j)\psi }d\psi. 
\end{equation}
Hence the matrix elements of the unitary irreducible representation of the
Euclidean group $E(2)$ with weight $p$ can be obtained from those of $SU(2)$
by the following contraction formula 
\begin{equation}
\label{eq: c2}t_{kj}^p (\phi , \rho , \zeta )= \lim _{l\rightarrow \infty
}t_{kj}^l(\phi ,p\rho /l,\zeta-\phi ). 
\end{equation}


\newpage

\begin{thebibliography}{9}
\bibitem{1}  Ahmedov, H and Duru, I. H., {\it J. Phys. A: Math. Gen., }, $%
{\bf 31}$, 5742 (1998).

\bibitem{2}  Faddev, L. D. and Takhtajan, L. A., {\it Lect. Notes Phys.}, $%
{\bf 246}$, 166 (1986).

\bibitem{3}  Vaksman, L. L. and Soibelman, Y. S., {\it Funkt. Anal. i
Prilozhen.}, ${\bf 22}$, 1 (1988).

\bibitem{5}  Vaksman, L. L., Korogodski, L. I., {\it Dokl. Akad. Nauk 
SSSR}, {\bf 304}, 
1036 (1989); Vaksman, L. L., {\it Dokl. Akad. Nauk SSSR}, {\bf 306}, 
269 (1989); Bonechi, F., Ciccoli, N., Giachetti, R., Sorace, E., Tarlini,
M., {\it Commun. Math. Phys. }, {\bf 175}, 161, (1996).

\bibitem{6}  In\"on\"u, E. and Wigner, E. P., {\it Proc. Nat. Akad. Sci.}, 
{\bf 39}, 510 (1956).

\bibitem{vil}  Vilenkin, N. Ja. and Klimyk, A. O., {\it Representation of
Lie Groups and Special Functions}, {\bf 3}, Dordrecht: Kluwer Akad. Publ.,
1992.

\bibitem{cel}  Celeghini, E., Giachetti, R., Sorace, E., Tarlini, M., {\it %
J. Phys. A : Math. Gen. }, {\bf 21}, 2548 (1990).

\bibitem{dob}  Dobrowski, L. and Sobczyk, J, {\it Lett. Math. Phys. }, {\bf %
32}, 249 (1994).

\bibitem{9}  Vilenkin, N. Ja., {\it Special Functions and Theory of Group
Representations}, Translation of Math. Monogr. {\bf 22}, Amer. Math. Soc.
Providence, Rhode Island.
\end{thebibliography}
\end{document}